\documentclass{spie_fickus}

\usepackage{amsthm}
\usepackage{amssymb}
\usepackage{amsmath}
\usepackage{graphicx}

\newcommand{\bbC}{\mathbb{C}}
\newcommand{\bbF}{\mathbb{F}}
\newcommand{\bbH}{\mathbb{H}}
\newcommand{\bbR}{\mathbb{R}}

\newcommand{\bfA}{\mathbf{A}}

\newcommand{\bfB}{\mathbf{B}}

\newcommand{\bfI}{\mathbf{I}}
\newcommand{\bfJ}{\mathbf{J}}

\newcommand{\bfP}{\mathbf{P}}
\newcommand{\bfQ}{\mathbf{Q}}

\newcommand{\bfU}{\mathbf{U}}
\newcommand{\bfV}{\mathbf{V}}

\newcommand{\bfx}{\mathbf{x}}

\newcommand{\bfy}{\mathbf{y}}

\newcommand{\bfone}{\boldsymbol{1}}
\newcommand{\bfzero}{\boldsymbol{0}}

\newcommand{\bfphi}{\boldsymbol{\varphi}}

\newcommand{\bfPhi}{\boldsymbol{\Phi}}

\newcommand{\bfSigma}{\boldsymbol{\Sigma}}

\newcommand{\calU}{\mathcal{U}}

\newcommand{\rmc}{\mathrm{c}}

\newcommand{\dist}{\operatorname{dist}}
\newcommand{\rank}{\operatorname{rank}}
\newcommand{\Span}{\operatorname{span}}
\newcommand{\Tr}{\operatorname{Tr}}

\newcommand{\Fro}{\mathrm{Fro}}

\newcommand{\abs}[1]{|{#1}|}

\newcommand{\biggparen}[1]{\biggl({#1}\biggr)}

\newcommand{\bigbracket}[1]{\bigl[{#1}\bigr]}

\newcommand{\set}[1]{\{{#1}\}}

\newcommand{\norm}[1]{\|{#1}\|}

\newcommand{\biggnorm}[1]{\biggl\|{#1}\biggr\|}

\newcommand{\ip}[2]{\langle{#1},{#2}\rangle}

\newtheorem{theorem}{Theorem}[section]

\newtheorem{lemma}[theorem]{Lemma}

\theoremstyle{definition}

\newtheorem{definition}[theorem]{Definition}

\newcommand{\romani}{\renewcommand{\labelenumi}{(\roman{enumi})}}

\title{A brief introduction to equi-chordal\\and equi-isoclinic tight fusion frames}

\author{Matthew Fickus\supit{a}, John Jasper\supit{b}, Dustin~G.~Mixon\supit{a}, Cody E.~Watson\supit{a}
\skiplinehalf
\supit{a}
Department of Mathematics and Statistics, Air Force Institute of Technology\\Wright-Patterson Air Force Base, Ohio 45433\\
\supit{b}
Department of Mathematical Sciences, University of Cincinnati, Cincinnati, OH 45221
}
\authorinfo{Send correspondence to Matthew Fickus: E-mail: Matthew.Fickus@gmail.com}

\begin{document}
\maketitle

\begin{abstract}
Equi-chordal and equi-isoclinic tight fusion frames (ECTFFs and EITFFs) are both types of optimal packings of subspaces in Euclidean spaces.
In the special case where these subspaces are one-dimensional, ECTFFs and EITFFs both correspond to types of optimal packings of lines known as equiangular tight frames.
In this brief note, we review some of the fundamental ideas and results concerning ECTFFs and EITFFs.
\end{abstract}

\section{Introduction}

In various settings, the following problem arises:
given a $d$-dimensional Hilbert space $\bbH$, and positive integers $c$ and $n$ with $c\leq d$,
how should we choose $n$ subspaces $\set{\calU_j}_{j=1}^{n}$ of $\bbH$, each of dimension $c$, such that the minimum pairwise distance between these subspaces is as large as possible?
That is, what are the optimal ways of packing $n$ points in the \textit{Grassmannian space} that consists of all $c$-dimensional subspaces of $\bbH$?
This is the central problem of a seminal paper by Conway, Hardin and Sloane~\cite{ConwayHS96}.

Of course, the answer here depends on how we define the distance between two such subspaces.
One popular choice is to let $\bfP_j:\bbH\rightarrow\bbH$ be the orthogonal projection operator onto $\calU_j$, and to define the \textit{chordal distance} between $\calU_j$ and $\calU_{j'}$ to be $\frac1{\sqrt{2}}$ times the Frobenius (Hilbert-Schmidt) distance between their projections.
In that setting, our goal is to maximize
\begin{equation}
\label{equation.chordal distance}
\dist_{\rmc}^2(\set{\calU_j}_{j=1}^{n}):=\tfrac1{2}\min_{j\neq j'}\norm{\bfP_j-\bfP_{j'}}_{\Fro}^2.
\end{equation}
Much previous work on this problem focuses on the special case where $c=1$.
In this case, letting $\bfphi_j$ be any unit vector in the $1$-dimensional subspace $\calU_j$,
we have $\bfP_j=\bfphi_j^{}\bfphi_j^{*}$ and, as detailed in a later section,
\eqref{equation.chordal distance} becomes
\begin{equation*}
\dist_{\rmc}^2(\set{\calU_j}_{j=1}^{n})
=\tfrac1{2}\min_{j\neq j'}\norm{\bfphi_j^{}\bfphi_j^{*}-\bfphi_{j'}^{}\bfphi_{j'}^{*}}_{\Fro}^2
=1-\max_{j\neq j'}\abs{\ip{\bfphi_j}{\bfphi_{j'}}}^2.
\end{equation*}
That is, our goal in this case is to choose unit vectors $\set{\bfphi_j}_{j=1}^{n}$ in $\bbH$ of minimal \textit{coherence}
$\max_{j\neq j'}\abs{\ip{\bfphi_j}{\bfphi_{j'}}}$.
Here, the \textit{Welch bound}~\cite{Welch74} states that for any positive integers $n\geq d$, any unit vectors \smash{$\set{\bfphi_j}_{j=1}^{n}$} in $\bbH$ satisfy
\begin{equation}
\label{equation.Welch bound}
\max_{j\neq j'}\abs{\ip{\bfphi_j}{\bfphi_{j'}}}^2
\geq\tfrac{(n-d)}{d(n-1)}.
\end{equation}
Moreover, it is well-known~\cite{StrohmerH03} that \smash{$\set{\bfphi_j}_{j=1}^{n}$} achieves equality in this bound if and only if it is an \textit{equiangular tight frame} (ETF) for $\bbH$,
namely when the value of $\abs{\ip{\bfphi_j}{\bfphi_{j'}}}$ is constant over all $j\neq j'$ (equiangularity) and there exists $\alpha>0$ such that $\sum_{j=1}^{n}\abs{\ip{\bfphi_j}{\bfx}}^2=\alpha\norm{\bfx}^2$ for all $\bfx\in\bbH$ (tightness).

In this short paper, we review some known results about these concepts,
including how the Welch bound~\eqref{equation.Welch bound} generalizes in the $c>1$ case to the \textit{simplex bound} of Conway, Hardin and Sloane~\cite{ConwayHS96}, which itself arises from a classical result by Rankin~\cite{Rankin55}.
Nearly all of these ideas are taken from Conway, Hardin and Sloane~\cite{ConwayHS96}, or from more recent investigations into these same ideas,
notably an article by Dhillon, Heath, Strohmer and Tropp~\cite{DhillonHST08},
as well as one by Kutyniok, Pezeshki, Calderbank and Liu~\cite{KutyniokPCL09}.

As we shall see, in order to achieve equality in the simplex bound, a sequence $\set{\calU_j}_{j=1}^{n}$ of $c$-dimensional subspaces of $\bbH$ is necessarily a \textit{tight fusion frame}~\cite{CasazzaK04}, that is, there necessarily exists $\alpha>0$ such that $\sum_{j=1}^{n}\bfP_j=\alpha\bfI$.
In particular, $\set{\calU_j}_{j=1}^{n}$ achieves equality in the simplex bound if and only if it is an \textit{equi-chordal tight fusion frame} (ECTFF) for $\bbH$~\cite{KutyniokPCL09}.
Refining that analysis leads to \textit{equi-isoclinic} subspaces~\cite{LemmensS73b},
and in particular, to \textit{equi-isoclinic tight fusion frames} (EITFFs).
In the special case when $c=1$,
both ECTFFs and EITFFs reduce to ETFs.

In the next section, we introduce the notation and basic concepts that we will need.
In Section~3, we review some classical results of Rankin~\cite{Rankin55} concerning optimal packings of points on real unit spheres.
Conway, Hardin and Sloane used these results to prove their simplex bound and \textit{orthoplex bound}, and we review their approach in the fourth section.
In Section~5, we provide a streamlined proof of the simplex bound,
and use it to explain why ECTFFs and EITFFs are the optimal solutions to two variants of the subspace packing problem.
In the final section, we interpret ECTFFs and EITFFs in terms of the \textit{principal angles} between two subspaces.

\section{Preliminaries}

Let $\bbH$ be a $d$-dimensional Hilbert space over a field $\bbF$ which is either $\bbR$ or $\bbC$.
The \textit{synthesis operator} of a finite sequence of vectors $\set{\bfphi_j}_{j=1}^{n}$ in $\bbH$ is the operator $\bfPhi:\bbF^n\rightarrow\bbH$, $\bfPhi\bfy:=\sum_{j=1}^{n}\bfy(n)\bfphi_j$.
Its adjoint is the \textit{analysis operator} $\bfPhi^*:\bbH\rightarrow\bbF^n$, $(\bfPhi^*\bfx)=\ip{\bfphi_j}{\bfx}$ where the inner product here is taken to be conjugate-linear in its first argument.
In the special case where $\bbH=\bbF^d$, which here is always assumed to be equipped with the standard (complex) dot product, $\bfPhi$ is the $d\times n$ matrix whose $j$th column is $\bfphi_j$, and $\bfPhi^*$ is its conjugate-transpose.
Composing these two operators gives the corresponding \textit{frame operator} $\bfPhi\bfPhi^*:\bbH\rightarrow\bbH$,
$\bfPhi\bfPhi^*\bfx=\sum_{j=1}^{n}\ip{\bfphi_j}{\bfx}\bfphi_j$,
and \textit{Gram matrix} $\bfPhi^*\bfPhi:\bbF^n\rightarrow\bbF^n$,
an $n\times n$ matrix whose $(j,j')$th entry is $(\bfPhi^*\bfPhi)(j,j')=\ip{\bfphi_j}{\bfphi_{j'}}$ for all $j,j'$.

A sequence of unit vectors $\set{\bfphi_j}_{j=1}^{n}$ in $\bbH$ is \textit{equiangular} if there exists $\beta\geq0$ such that $\abs{\ip{\bfphi_j}{\bfphi_{j'}}}^2=\beta$ for all $j\neq j'$;
to clarify, under this definition, it is not the vectors themselves that are necessarily equiangular,
but rather, the lines they span.
Meanwhile, $\set{\bfphi_j}_{j=1}^{n}$ is a tight frame for $\bbH$ if there exists $\alpha>0$ such that $\alpha\norm{\bfx}^2=\sum_{j=1}^{n}\abs{\ip{\bfphi_j}{\bfx}}^2$ for all $\bfx\in\bbH$.
This requires that $\set{\bfphi_j}_{j=1}^{n}$ spans $\bbH$.
In particular, we need $n\geq\dim(\bbH)$.
By the polarization identity, $\set{\bfphi_j}_{j=1}^{n}$ is tight with tight frame constant $\alpha>0$ precisely when
\begin{equation*}
\alpha\bfI
=\bfPhi\bfPhi^*
=\sum_{j=1}^{n}\bfphi_j^{}\bfphi_j^*.
\end{equation*}
Here, $\bfphi_j^*:\bbH\rightarrow\bbF$ is the linear functional $\bfphi_j^*\bfx:=\ip{\bfphi_j}{\bfx}$ that arises as the adjoint of viewing $\bfphi_j$ as an operator $\bfphi_j:\bbF\rightarrow\bbH$, $\bfphi_j(c):=c\bfphi_j$.
In the special case where $\bfphi_j$ is unit norm,
the operator $\bfphi_j^{}\bfphi_j^*$ is the orthogonal projection operator onto the line $\Span\set{\bfphi_j}$.
That is, a \textit{unit norm tight frame} corresponds to a collection of rank-one orthogonal projection operators which sum to a scalar multiple of the identity operator.
In this case, the tight frame constant $\alpha$ is necessarily the \textit{redundancy} $\frac nd$ of the frame since the diagonal entries of the Gram matrix $\bfPhi^*\bfPhi$ are all $1$ and so $\alpha d=\Tr(\alpha\bfI)=\Tr(\bfPhi\bfPhi^*)=\Tr(\bfPhi^*\bfPhi)=n$.

When $\set{\bfphi_j}_{j=1}^{n}$ is both equiangular and a tight frame for $\bbH$ we say it is an \textit{equiangular tight frame} (ETF) for $\bbH$.
As detailed in a later section, being an ETF is equivalent to achieving equality in the Welch bound~\cite{Welch74}.

The theory of unit norm tight frames naturally generalizes to subspaces $\set{\calU_j}_{j=1}^{n}$ of $\bbH$ of dimension $c>1$.
Here, for each $j=1,\dotsc,n$, let $\bfPhi_j:\bbF^c\rightarrow\bbH$ be the synthesis operator for an orthonormal basis $\set{\bfphi_{j,k}}_{k=1}^{c}$ of $\calU_j$.
For example, when $\bbH=\bbF_d$, $\bfPhi_j$ is a $d\times c$ matrix whose $k$th column is $\bfphi_{j,k}$.
The fact that $\set{\bfphi_{j,k}}_{k=1}^{c}$ is an orthonormal basis for $\calU_j$ implies that $\bfPhi_j^*\bfPhi_j^{}=\bfI$ and that $\bfP_j=\bfPhi_j^{}\bfPhi_j^{*}$ is the orthogonal projection operator onto $\calU_j$.
The Gram matrix $\bfPhi^*\bfPhi$ of the $nc$ concatenated vectors $\set{\bfphi_{j,k}}_{j=1}^{n}\,_{k=1}^{c}$ is called a \textit{fusion Gram matrix} of $\set{\calU_j}_{j=1}^{n}$,
and can be regarded as an $n\times n$ array of $c\times c$ submatrices.
Specifically, for any $j,j'=1,\dotsc,n$, the $(j,j')$th block of $\bfPhi^*\bfPhi$ is the \textit{cross-Gramian} $\bfPhi_j^*\bfPhi_{j'}^{}$ whose $(k,k')$th entry is $(\bfPhi_j^*\bfPhi_{j'}^{})(k,k')=\ip{\bfphi_{j,k}}{\bfphi_{j',k'}}$.
Meanwhile, the frame operator $\bfPhi\bfPhi^*$ of $\set{\bfphi_{j,k}}_{j=1}^{n}\,_{k=1}^{c}$ is called the \textit{fusion frame operator} of $\set{\calU_j}_{j=1}^{n}$,
and is the sum of the corresponding rank-$c$ orthogonal projection operators $\set{\bfP_j}_{j=1}^{n}$,
\begin{equation}
\label{equation.fusion frame operator}
\bfPhi\bfPhi^*
=\sum_{j=1}^{n}\sum_{k=1}^{c}\bfphi_{j,k}^{}\bfphi_{j,k}^{*}
=\sum_{j=1}^{n}\bfPhi_j^{}\bfPhi_j^{*}
=\sum_{j=1}^{n}\bfP_j.
\end{equation}

Note here that the fusion Gram matrix of a fusion frame is not unique, as it depends on the particular choice of orthonormal basis for each $\calU_j$.
(To be clear, the diagonal blocks of any fusion Gram matrix are all $\bfI$ regardless.)
However, the fusion frame operator is unique, as the orthogonal projection operator onto $\calU_j$ is invariant with respect to such choices.
More precisely, each $\bfPhi_j$ is unique up to right-multiplication by unitary $c\times c$ matrices, meaning $\bfPhi_j^{}\bfPhi_j^*$ is unique while $\bfPhi_j^*\bfPhi_{j'}^{}$ is only unique up to right- and left-multiplication by unitaries.
Because of this, it is more natural to generalize the notion of a tight frame to this ``fusion" setting than it is to generalize the notion of an equiangular frame.
Indeed, as noted in the introduction, $\set{\calU_j}_{j=1}^{n}$ is a \textit{tight fusion frame} for $\bbH$ if the corresponding fusion frame operator~\eqref{equation.fusion frame operator} is $\alpha\bfI$ for some $\alpha>0$; in this case taking the trace of this equation gives
$\alpha d=\Tr(\alpha\bfI)=\sum_{j=1}^{n}\Tr(\bfP_j)=nc$ and so $\alpha$ is necessarily $\tfrac{nc}d$.
To generalize equiangularity, we want some way to take the ``modulus" of the off-diagonal cross-Gramians $\bfPhi_j^*\bfPhi_{j'}$ that is invariant with respect to left- and right-multiplication by unitaries.
As we shall see, there is more than one option here:
taking the Frobenius norms of these matrices leads to the notion of \textit{equi-chordal} subspaces, while taking the induced $2$-norms of these matrices leads to \textit{equi-isoclinic} subspaces.

\section{Optimal packings on the sphere}

Conway, Hardin and Sloane~\cite{ConwayHS96} give two distinct sufficient conditions for a sequence $\set{\calU_j}_{j=1}^{n}$ of $c$-dimensional subspaces  of a $d$-dimensional Hilbert space $\bbH$ to form an optimal packing with respect to the chordal distance~\eqref{equation.chordal distance}.
These results are often referred to as the \textit{simplex bound} and \textit{orthoplex bound}.
The traditional proofs of these results depend on Rankin's earlier work concerning optimal packings on real spheres~\cite{Rankin55},
that is, ways to arrange $n$ unit vectors in $\bbR^d$ whose minimum pairwise distance is as large as possible.
In particular, the traditional proof of the simplex bound depends on the following concept:

\begin{definition}
A sequence of unit vectors $\set{\bfphi_j}_{j=1}^{n}$ in a Hilbert space $\bbH$ is a called a \textit{regular simplex} if there exists some integer $n\geq 2$ such that $\ip{\bfphi_j}{\bfphi_{j'}}=-\frac1{n-1}$ for all $j\neq j'$.
\end{definition}

That is, $\set{\bfphi_j}_{j=1}^{n}$ is a regular simplex if and only if its Gram matrix is
\smash{$\bfPhi^*\bfPhi=\frac{n}{n-1}\bfI-\frac1{n-1}\bfJ$}.
The eigenvalues of such a Gram matrix are $0$ and $\frac{n}{n-1}$ with eigenspaces $\Span\set{\bfone}$ and \smash{$\bfone^\perp$}, respectively.
In particular, the null space of $\bfPhi$ is
$\ker(\bfPhi)=\ker(\bfPhi^*\bfPhi)=\Span\set{\bfone}$,
meaning \smash{$\sum_{j=1}^{n}\bfphi_j=\bfPhi\bfone=\bfzero$}.
Moreover, the dimension of $\Span\set{\bfphi_j}_{j=1}^{n}$ is $\rank(\bfPhi)=n-1$.

Regular simplices are optimal packings on real unit spheres.
For example, in $\bbR^3$, an optimal packing of two unit vectors consists of a pair of antipodal unit vectors, whereas an optimal packing of three unit vectors consists of three equally-spaced vectors in a great circle, and an optimal packing of four unit vectors forms a tetrahedron, namely regular simplices with $n=2$, $3$ and $4$ vectors, respectively.
To prove this formally, note that for any finite sequence $\set{\bfphi_j}_{j=1}^{n}$ of unit norm vectors in a real Hilbert space,
\begin{equation*}
0
\leq\biggnorm{\sum_{j=1}^{n}\bfphi_j}^2
=\sum_{j=1}^{n}\sum_{j'=1}^{n}\ip{\bfphi_j}{\bfphi_{j'}}
=n+\sum_{j=1}^{n}\sum_{\substack{j'=1\\j'\neq j}}^{n}\ip{\bfphi_j}{\bfphi_{j'}}
\leq n+n(n-1)\max_{j\neq j'}\ip{\bfphi_j}{\bfphi_{j'}}.
\end{equation*}
That is, for any such vectors we have
$-\frac1{n-1}\leq\max_{j\neq j'}\ip{\bfphi_j}{\bfphi_{j'}}$.
Moreover, both of the above inequalities hold with equality if and only if $\ip{\bfphi_j}{\bfphi_{j'}}=-\frac1{n-1}$, that is, if and only if $\set{\bfphi_j}_{j=1}^{n}$ is a regular simplex.
To relate this to optimal packings, note that for any unit vectors $\set{\bfphi_j}_{j=1}^{n}$ in a real Hilbert space,
\begin{equation}
\label{equation.packing in terms of coherence}
\min_{j\neq j'}\norm{\bfphi_j-\bfphi_{j'}}^2
=\min_{j\neq j'}2(1-\ip{\bfphi_j}{\bfphi_{j'}})
=2(1-\max_{j\neq j'}\ip{\bfphi_j}{\bfphi_{j'}}).
\end{equation}
We summarize these facts in the following result:

\begin{lemma}
\label{lemma.Rankin simplex bound}
For any finite sequence $\set{\bfphi_j}_{j=1}^{n}$ of unit norm vectors in a real Hilbert space,
\begin{equation*}
-\tfrac1{n-1}
\leq\max_{j\neq j'}\ip{\bfphi_j}{\bfphi_{j'}}
=1-\tfrac12\min_{j\neq j'}\norm{\bfphi_j-\bfphi_{j'}}^2,
\end{equation*}
where equality holds if and only if $\set{\bfphi_j}_{j=1}^{n}$ is a regular simplex.
\end{lemma}

We need $n\leq d+1$ in order to achieve equality here:
if $\set{\bfphi_j}_{j=1}^{n}$ is a regular simplex that lies in $\bbH$ then
$n-1=\Span(\bfPhi)\leq\dim(\bbH)=d$.
For example, the optimal packing of $n=5$ unit vectors in $\bbR^3$ cannot be a regular simplex.
In fact, as we now discuss, when $n>d+1$ one can prove a bound that is stronger than that given in Lemma~\ref{lemma.Rankin simplex bound}.

In particular, when $n>d+1$, we have \smash{$\max_{j\neq j'}\ip{\bfphi_j}{\bfphi_{j'}}\geq0$} for any unit vectors $\set{\bfphi_j}_{j=1}^{n}$ in a real $d$-dimensional Hilbert space $\bbH$.
For an elegant proof of this fact~\cite{Chapman10},
note that $\dim(\ker(\bfPhi))\geq 2$ when $n\geq d+2$,
implying there exist nontrivial, nonnegative vectors $\bfy_1,\bfy_2\in\ker(\bfPhi)$ with disjoint support;
thus,
\begin{equation*}
0
=\ip{\bfzero}{\bfzero}
=\ip{\bfPhi\bfy_1}{\bfPhi\bfy_2}
=\sum_{j=1}^n\sum_{\substack{j'=1\\j'\neq j}}^n\bfy_1(j)\bfy_2(j')\ip{\bfphi_j}{\bfphi_{j'}},
\end{equation*}
where $\bfy_1(j)\bfy_2(j')$ is nonnegative for all $j,j'$, and is strictly positive for at least one pair $j\neq j'$;
as such we cannot have $\ip{\bfphi_j}{\bfphi_{j'}}<0$ for all $j\neq j'$.
In summary, when combined with~\eqref{equation.packing in terms of coherence},
we have the following result:

\begin{lemma}
\label{lemma.Rankin orthoplex bound}
For any positive integers $n$ and $d$ with $n\geq d+2$,
and any finite sequence $\set{\bfphi_j}_{j=1}^{n}$ in a real $d$-dimensional Hilbert space,
\begin{equation*}
0\leq\max_{j\neq j'}\ip{\bfphi_j}{\bfphi_{j'}}
=1-\tfrac12\min_{j\neq j'}\norm{\bfphi_j-\bfphi_{j'}}^2.
\end{equation*}
\end{lemma}

As noted by Rankin~\cite{Rankin55}, equality can be achieved here with $2d$ vectors by choosing $\set{\bfphi_j}_{j=1}^n$ to be an orthonormal basis along with its antipodes, i.e., $\set{\pm\delta_j}_{j=1}^{d}$;
such a sequence of vectors is known as an \textit{orthoplex}.

\section{The simplex and orthoplex bounds}

To obtain the simplex and orthoplex bounds of Conway, Hardin and Sloane,
we apply Lemmas~\ref{lemma.Rankin simplex bound} and~\ref{lemma.Rankin orthoplex bound} to the normalized traceless components of orthogonal projection operators.
To be precise, let $\set{\calU_j}_{j=1}^{n}$ be $c$-dimensional subspaces of $\bbF^n$ where $\bbF$ is either $\bbR$ or $\bbC$, and for each $j=1,\dotsc,n$, let $\bfP_j$ be the $n\times n$ matrix which is the orthogonal projection operator onto $\calU_j$.
The operators $\set{\bfP_j}_{j=1}^{n}$ lie in the real Hilbert space of all $d\times d$ self-adjoint matrices; here, the inner product is the Frobenius (Hilbert-Schmidt) inner product $\ip{\bfA}{\bfB}_{\Fro}:=\Tr(\bfA^*\bfB)$, and this space has dimension \smash{$\frac{d(d+1)}{2}$} or $d^2$ depending on whether $\bbF$ is $\bbR$ or $\bbC$, respectively.
A $d\times d$ self-adjoint matrix is \textit{traceless} if its trace is zero,
namely if it lies in the orthogonal complement of $\bfI$.
Since $\Tr(\bfP_j)=\dim(\calU_j)=c$ for all $j$, the traceless component of any $\bfP_j$ is
\begin{equation*}
\bfP_j-\tfrac{\ip{\bfI}{\bfP_j}_{\Fro}}{\ip{\bfI}{\bfI}_{\Fro}}\bfI
=\bfP_j-\tfrac{\Tr(\bfP_j)}{\Tr(\bfI)}\bfI
=\bfP_j-\tfrac{c}{d}\bfI.
\end{equation*}
To continue, note
\begin{equation}
\label{equation.deriving inner product of traceless}
\ip{\bfP_j-\tfrac{c}{d}\bfI}{\bfP_{j'}-\tfrac{c}{d}\bfI}_{\Fro}
=\Tr[(\bfP_j-\tfrac{c}{d}\bfI)(\bfP_{j'}-\tfrac{c}{d}\bfI)]
=\Tr(\bfP_j\bfP_{j'})-2\tfrac{c^2}{d}+\tfrac{c^2}{d}
=\ip{\bfP_j}{\bfP_{j'}}_{\Fro}-\tfrac{c^2}{d}
\end{equation}
for all $j,j'$.
In particular, $\norm{\bfP_j-\tfrac{c}{d}\bfI}_{\Fro}^2=c-\tfrac{c^2}{d}=\tfrac{c(d-c)}{d}$ for all $j$, meaning
\begin{equation*}
\set{\bfQ_j}_{j=1}^{n},
\quad
\bfQ_j:=\bigbracket{\tfrac{d}{c(d-c)}}^{\frac12}(\bfP_j-\tfrac{c}{d}\bfI),
\end{equation*}
is a normalized sequence of vectors in a real Hilbert space of dimension $\frac{d(d+1)}{2}-1$ or $d^2-1$, depending on whether $\bbF$ is $\bbR$ or $\bbC$, respectively.
Moreover, by~\eqref{equation.deriving inner product of traceless},
\begin{equation*}
\ip{\bfQ_j}{\bfQ_{j'}}_{\Fro}
=\tfrac{d}{c(d-c)}(\ip{\bfP_j}{\bfP_{j'}}_{\Fro}-\tfrac{c^2}{d}),
\end{equation*}
for all $j,j'$.
As such, applying Lemma~\ref{lemma.Rankin simplex bound} to $\set{\bfQ_j}_{j=1}^{n}$ and recalling~\eqref{equation.chordal distance} then gives
\begin{equation}
\label{equation.proving simplex from Rankin}
-\tfrac1{n-1}
\leq\tfrac{d}{c(d-c)}(\max_{j\neq j'}\ip{\bfP_j}{\bfP_{j'}}_{\Fro}-\tfrac{c^2}{d})
=1-\tfrac12\tfrac{d}{c(d-c)}\min_{j\neq j'}\norm{\bfP_j-\bfP_{j'}}_{\Fro}^2
=1-\tfrac{d}{c(d-c)}\dist_{\rmc}^2(\set{\calU_j}_{j=1}^{n}).
\end{equation}
Rearranging this expression gives the simplex bound of Conway, Hardin and Sloane~\cite{ConwayHS96}:
\begin{equation}
\label{equation.simplex bound}
\dist_{\rmc}^2(\set{\calU_j}_{j=1}^{n})
\leq \tfrac{c(d-c)}{d}\tfrac{n}{n-1}.
\end{equation}
If we instead solve for $\max_{j\neq j'}\ip{\bfP_j}{\bfP_{j'}}_{\Fro}$ above, we obtain a lower bound on that quantity which is equivalent to the simplex bound~\eqref{equation.simplex bound}:
\begin{equation}
\label{equation.Welch-simplex bound}
\max_{j\neq j'}\ip{\bfP_j}{\bfP_{j'}}_{\Fro}
\geq\tfrac{c^2}{d}-\tfrac1{n-1}\tfrac{c(d-c)}{d}
=\tfrac{c}{d(n-1)}[(n-1)c-(d-c)]
=\tfrac{c(nc-d)}{d(n-1)}.
\end{equation}
Note that by Lemma~\ref{lemma.Rankin simplex bound},
\eqref{equation.simplex bound} and \eqref{equation.Welch-simplex bound} hold precisely when $\set{\bfQ_j}_{j=1}^{n}$ forms a simplex in the space of all (traceless) self-adjoint operators,
namely when $\ip{\bfP_j}{\bfP_{j'}}_{\Fro}=\tfrac{c(nc-d)}{d(n-1)}$ for all $j\neq j'$.
Recall this can only happen if
\begin{equation*}
\bfzero
=\sum_{j=1}^{n}\bfQ_j
=\bigbracket{\tfrac{d}{c(d-c)}}^{\frac12}\sum_{j=1}^{n}(\bfP_j-\tfrac{c}{d}\bfI)
=\bigbracket{\tfrac{d}{c(d-c)}}^{\frac12}\biggparen{\sum_{j=1}^{n}\bfP_j-\tfrac{nc}{d}\bfI},
\end{equation*}
namely only if $\set{\calU_j}_{j=1}^{n}$ is a tight fusion frame for $\bbF^d$.

An important fact seemingly overlooked by Conway, Hardin and Sloane is that the simplex bound~\eqref{equation.simplex bound}, when equivalently reexpressed as~\eqref{equation.Welch-simplex bound}, is a generalization of the Welch bound~\eqref{equation.Welch bound}.
To see this, note that when $c=1$, then for each $j=1,\dotsc,n$ we have $\bfP_j=\bfphi_j^{}\bfphi_j^{}$ where $\bfphi_j$ is a unit vector in the $1$-dimensional subspace $\calU_j$.
Moreover, by cycling a trace (and realizing that the trace of a scalar is itself), we have
\begin{equation*}
\ip{\bfP_j}{\bfP_{j'}}_{\Fro}
=\Tr(\bfP_j\bfP_{j'})
=\Tr(\bfphi_j^{}\bfphi_j^{*}\bfphi_{j'}^{}\bfphi_{j'}^{*})
=\Tr(\bfphi_{j'}^{*}\bfphi_j^{}\bfphi_j^{*}\bfphi_{j'}^{})
=\Tr(\ip{\bfphi_{j'}}{\bfphi_{j}}\ip{\bfphi_{j}}{\bfphi_{j'}})
=\abs{\ip{\bfphi_{j}}{\bfphi_{j'}}}^2
\end{equation*}
for all $j,j'$.
As such, in the $c=1$ case, \eqref{equation.Welch-simplex bound} reduces to
\smash{$\max_{j\neq j'}\abs{\ip{\bfphi_{j}}{\bfphi_{j'}}}^2\geq\tfrac{n-d}{d(n-1)}$}, namely \eqref{equation.Welch bound}.
This realization inspired the discussion given in the next section;
there, we obtain a more direct proof of~\eqref{equation.Welch-simplex bound} by generalizing a modern proof of the Welch bound.

To conclude this section, we note that when $\bbF=\bbR$ and \smash{$n>\frac{d(d+1)}{2}$},
the $n$ vectors $\set{\bfQ_j}_{j=1}^{n}$ cannot form a simplex in the \smash{$[\frac{d(d+1)}{2}-1]$}-dimensional space of $d\times d$ real symmetric traceless operators.
In the special case where $c=1$, this fact is closely related to the \textit{Gerzon bound}~\cite{LemmensS73}, which states that the maximum number of equiangular lines in $\bbR^d$ is \smash{$\frac{d(d+1)}{2}$}.
In this regime,
we can instead apply Lemma~\ref{lemma.Rankin orthoplex bound} to $\set{\bfQ_j}_{j=1}^{n}$ to obtain the following alternative to \eqref{equation.proving simplex from Rankin}:
\begin{equation*}
0
\leq\tfrac{d}{c(d-c)}(\max_{j\neq j'}\ip{\bfP_j}{\bfP_{j'}}_{\Fro}-\tfrac{c^2}{d})
=1-\tfrac{d}{c(d-c)}\dist_{\rmc}^2(\set{\calU_j}_{j=1}^{n}).
\end{equation*}
In the case where $\bbF=\bbC$, these same bounds apply when $n>d^2$.
When rearranged, this yields Conway, Hardin and Sloane's orthoplex bound~\cite{ConwayHS96}:
$\dist_{\rmc}^2(\set{\calU_j}_{j=1}^{n})\leq\tfrac{c(d-c)}{d}$.
Equivalently,
\smash{$\max_{j\neq j'}\ip{\bfP_j}{\bfP_{j'}}_{\Fro}\geq\tfrac{c^2}{d}$}.
In the special case where $c=1$ and $\bfP_j=\bfphi_j^{}\bfphi_j^{*}$ for all $j$,
this becomes
\begin{equation*}
\max_{j\neq j'}\abs{\ip{\bfphi_{j}}{\bfphi_{j'}}}^2\geq\tfrac1d,
\end{equation*}
a bound met by \textit{mutually unbiased bases} as well as by other interesting constructions, such as the union of a standard basis and a harmonic ETF arising from a Singer difference set~\cite{BodmannH16}.

\section{Equi-chordal and equi-isoclinic tight fusion frames}

In this section, we give an alternative derivation of~\eqref{equation.Welch-simplex bound},
an inequality that is equivalent to Conway, Hardin and Sloane's simplex bound~\eqref{equation.simplex bound}.
This is not a new proof \textit{per se}: it essentially combines the main ideas of the previous section with those of the proof of Lemma~\ref{lemma.Rankin simplex bound}, while eliminating some unnecessary technicalities.

Recall from Section~2, that a sequence $\set{\calU_j}_{j=1}^{n}$ of $c$-dimensional subspaces of $\bbF^d$ forms a tight fusion frame for $\bbF^d$ if there exists $\alpha>0$ such that $\alpha\bfI=\sum_{j=1}^{n}\bfPhi_j^{}\bfPhi_j^{*}$ where each $\bfPhi_j$ is a $d\times c$ synthesis operator of an orthonormal basis for $\calU_j$.
Further recall that in this case, $\alpha$ is necessarily $\frac{nc}{d}$.
As such, for any sequence $\set{\calU_j}_{j=1}^{n}$ of $c$-dimensional subspaces of $\bbF^d$,
\begin{align*}
0
&\leq\norm{\bfPhi\bfPhi^*-\tfrac{nc}{d}\bfI}_{\Fro}^2\\
&=\Tr[(\bfPhi\bfPhi^*-\tfrac{nc}{d}\bfI)^2]\\
&=\Tr[(\bfPhi\bfPhi^*)^2]-2\tfrac{nc}{d}\Tr(\bfPhi\bfPhi^*)+\tfrac{n^2c^2}{d^2}\Tr(\bfI)\\
&=\Tr[(\bfPhi^*\bfPhi)^2]-2\tfrac{nc}{d}\Tr(\bfPhi^*\bfPhi)+\tfrac{n^2c^2}{d}\\
&=\norm{\bfPhi^*\bfPhi}_{\Fro}^2-\tfrac{n^2c^2}{d}.
\end{align*}
To continue simplifying, we express the Frobenius norm (entrywise $2$-norm) of the fusion Gram matrix in terms of the Frobenius norms of the cross-Gramians $\set{\bfPhi_j^*\bfPhi_{j'}^{}}_{j,j'=1}^n$, recalling that $\bfPhi_j^*\bfPhi_j=\bfI$ for all $j=1,\dotsc,n$:
\begin{equation}
\label{equation.direct derivation of generalized Welch}
0
\leq\norm{\bfPhi\bfPhi^*-\tfrac{nc}{d}\bfI}_{\Fro}^2
=\sum_{j=1}^{n}\sum_{j'=1}^{n}\norm{\bfPhi_j^*\bfPhi_{j'}^{}}_{\Fro}^2-\tfrac{n^2c^2}{d}
\leq n(n-1)\max_{j\neq j'}\norm{\bfPhi_j^*\bfPhi_{j'}^{}}_{\Fro}^2+nc-\tfrac{n^2c^2}{d}.
\end{equation}
Rearranging this inequality yields~\eqref{equation.Welch-simplex bound},
namely the inequality that is equivalent to the simplex bound.
Moreover, equality in~\eqref{equation.direct derivation of generalized Welch} only occurs when both inequalities hold with equality, namely when $\bfPhi\bfPhi^*=\tfrac{nc}{d}\bfI$ and $\norm{\bfPhi_j^*\bfPhi_{j'}^{}}_{\Fro}^2$ is constant over all $j\neq j'$.
Recall our first property here means $\set{\calU_j}_{j=1}^{n}$ is a tight fusion frame for $\bbF^d$.
Meanwhile, our second property is equivalent to $\set{\calU_j}_{j=1}^{n}$ being \textit{equi-chordal}, namely that the chordal distance between any two distinct subspaces is the same value:
\begin{align}
\dist_{\rmc}^2(\calU_j,\calU_{j'})
\nonumber
&=\tfrac12\norm{\bfP_j-\bfP_{j'}}_{\Fro}^2\\
\nonumber
&=\tfrac12\Tr[(\bfP_j-\bfP_{j'})^2]\\
\nonumber
&=\tfrac12[\Tr(\bfP_j)+\Tr(\bfP_j)-2\Tr(\bfP_j\bfP_{j'})]\\
\nonumber
&=c-\Tr(\bfPhi_j^{}\bfPhi_j^*\bfPhi_{j'}^{}\bfPhi_{j'}^*)\\
\nonumber
&=c-\Tr(\bfPhi_{j'}^*\bfPhi_j^{}\bfPhi_j^*\bfPhi_{j'}^{})\\
\label{equation.chordal distance between two subspaces}
&=c-\norm{\bfPhi_j^*\bfPhi_{j'}^{}}_{\Fro}^2.
\end{align}
That is, equality in~\eqref{equation.direct derivation of generalized Welch} is achieved precisely when $\set{\calU_j}_{j=1}^{n}$ is equi-chordal and is a tight fusion frame.
We summarize these facts in the following result:

\begin{theorem}
\label{theorem.equi-chordal}
Let $\set{\calU_j}_{j=1}^{n}$ be a sequence of $c$-dimensional subspaces of $\bbF^d$.
For each $j=1,\dotsc,n$, let $\bfPhi_j$ be the $d\times c$ synthesis operator of an orthonormal basis for $\calU_j$, and let $\bfP_j=\bfPhi_j^{}\bfPhi_j^*$ be the corresponding orthogonal projection operator.
Then
\begin{equation}
\label{equation.equi-chordal Welch bound}
\max_{j\neq j'}\norm{\bfPhi_j^*\bfPhi_{j'}^{}}_{\Fro}^2
\geq\tfrac{c(nc-d)}{d(n-1)},
\end{equation}
where equality holds if and only if $\set{\calU_j}_{j=1}^{n}$ is an equi-chordal tight fusion frame (ECTFF) for $\bbF^d$, namely if and only if both
\begin{enumerate}
\romani
\item
$\set{\calU_j}_{j=1}^{n}$ is a tight fusion frame,
namely there exists some $\alpha>0$ such that
$\displaystyle\sum_{j=1}^{n}\bfPhi_j^{}\bfPhi_j^{*}=\sum_{j=1}^{n}\bfP_j=\alpha\bfI$;
\item
$\set{\calU_j}_{j=1}^{n}$ is equi-chordal, namely there exists some $\beta\geq0$ such that $\Tr(\bfP_j\bfP_{j'})=\norm{\bfPhi_j^*\bfPhi_{j'}^{}}_{\Fro}^2=\beta$ for all $j\neq j'$,
or equivalently, that the squared chordal distance $\frac12\norm{\bfP_j-\bfP_{j'}}^2$ is constant over all $j\neq j'$.
\end{enumerate}
In this case, $\alpha$ and $\beta$ are necessarily $\frac{nc}{d}$ and $\tfrac{c(nc-d)}{d(n-1)}$, respectively.
\end{theorem}

We now further refine these ideas to obtain another bound that can only be achieved by \textit{equi-isoclinic} subspaces.
Recall that the Frobenius norm of a matrix is the $2$-norm of its singular values,
while its induced $2$-norm is the $\infty$-norm of its singular values, yielding the following bound:
writing $\bfA\in\bbC^{c\times c}$ as $\bfA=\bfU\bfSigma\bfV^*$ where $\bfU$ and $\bfV$ are unitary,
\begin{equation*}
\norm{\bfA}_{\Fro}^2
=\norm{\bfU\bfSigma\bfV^*}_{\Fro}^2
=\norm{\bfSigma}_{\Fro}^2
=\sum_{k=1}^{c}\sigma_k^2
\leq c\max\set{\sigma_k^2}_{k=1}^{c}
=c\norm{\bfA}_2^2.
\end{equation*}
Moreover, equality here is only achieved when the singular values $\set{\sigma_k}_{k=1}^{c}$ are constant, namely when there exists some $\sigma\geq0$ such that $\bfSigma=\sigma\bfI$.
As the singular values of $\bfA\in\bbC^{c\times c}$ are the square roots of the eigenvalues of $\bfA^*\bfA$ and $\bfA\bfA^*$, this occurs precisely when $\bfA^*\bfA=\sigma^2\bfI$, or equivalently, when $\bfA\bfA^*=\sigma^2\bfI$.
That is, for any $\bfA\in\bbC^{c\times c}$, $\norm{\bfA}_{\Fro}^2\leq c\norm{\bfA}_2^2$, and equality is achieved if and only if $\bfA$ is a scalar multiple of a unitary matrix.
Note that in this case, the value $\sigma$ is uniquely determined by the Frobenius norm of $\bfA$,
namely $\sigma^2=\frac1{c}\norm{\bfA}_{\Fro}^2$.

Using these ideas, we can continue~\eqref{equation.equi-chordal Welch bound} as
\begin{equation}
\label{equation.deriving EITFFs}
\tfrac{c(nc-d)}{d(n-1)}
\leq\max_{j\neq j'}\norm{\bfPhi_j^*\bfPhi_{j'}^{}}_{\Fro}^2
\leq c\max_{j\neq j'}\norm{\bfPhi_j^*\bfPhi_{j'}^{}}_{2}^2,
\end{equation}
obtaining the lower bound $\max_{j\neq j'}\norm{\bfPhi_j^*\bfPhi_{j'}^{}}_{2}^2
\geq\tfrac{nc-d}{d(n-1)}$.
Here, we have equality precisely when we have equality in~\eqref{equation.equi-chordal Welch bound} and $\bfPhi_j^*\bfPhi_{j'}^{}$ has constant singular values for any $j\neq j'$.
That is, equality in~\eqref{equation.deriving EITFFs} holds precisely when $\set{\calU_j}_{j=1}^{n}$ is an ECTFF and there exists $\sigma\geq0$ such that
$\bfPhi_{j'}^*\bfPhi_j^{}\bfPhi_j^*\bfPhi_{j'}^{}=\sigma^2\bfI$ for all $j\neq j'$;
note the value of $\sigma$ here is independent of $j,j'$, satisfying
$\sigma^2=\frac1{c}\norm{\bfPhi_j^*\bfPhi_{j'}^{}}_{\Fro}^2=\frac1c\beta=\tfrac{nc-d}{d(n-1)}$.
As we now discuss, sequences of subspaces with this special property are themselves a subject of interest.

In particular, a sequence of $c$-dimensional subspaces $\set{\calU_j}_{j=1}^{n}$ of $\bbF^d$ is called \textit{equi-isoclinic} if there exists $\sigma\geq0$ such that $\bfPhi_{j'}^*\bfPhi_j^{}\bfPhi_j^*\bfPhi_{j'}^{}=\sigma^2\bfI$ for all $j\neq j'$~\cite{LemmensS73b}.
Note that conjugating this expression by $\bfPhi_{j'}$ gives
\begin{equation}
\label{equation.equi-isoclinic projections}
\bfP_{j'}\bfP_{j}\bfP_{j'}
=\bfPhi_{j'}^{}\bfPhi_{j'}^*\bfPhi_j^{}\bfPhi_j^*\bfPhi_{j'}^{}\bfPhi_{j'}^{*}
=\sigma^2\bfPhi_{j'}^{}\bfPhi_{j'}^{*}
=\sigma^2\bfP_{j'}
\end{equation}
for all $j\neq j'$.
Conversely, as we now explain, if there exists some $\sigma\geq0$ such that the subspaces $\set{\calU_j}_{j=1}^{n}$ satisfy \eqref{equation.equi-isoclinic projections} for all $j\neq j'$,  then they are equi-isoclinic.
Indeed, \eqref{equation.equi-isoclinic projections} implies that
$\bfP_{j'}(\bfP_{j}\bfP_{j'}-\sigma^2\bfI)=\bfzero$,
meaning the range of $\bfP_{j}\bfP_{j'}-\sigma^2\bfI$ lies in $\ker(\bfP_{j'})=\ker(\bfPhi_{j'}^{}\bfPhi_{j'}^*)=\ker(\bfPhi_{j'}^*)$.
That is,
$\bfzero
=\bfPhi_{j'}^*(\bfP_{j}^{}\bfP_{j'}^{}-\sigma^2\bfI)
=(\bfPhi_{j'}^*\bfP_{j}^{}\bfPhi_{j'}^{}-\sigma^2\bfI)\bfPhi_{j'}^*$.
Since $\bfPhi_{j'}^*\bfPhi_{j'}^{}=\bfI$, multiplying this equation on the right by $\bfPhi_{j'}$ gives
$\bfPhi_{j'}^*\bfPhi_j^{}\bfPhi_j^*\bfPhi_{j'}^{}
=\bfPhi_{j'}^*\bfP_{j}^{}\bfPhi_{j'}^{}
=\sigma^2\bfI$.
We summarize these facts as follows:

\begin{theorem}
Let $\set{\calU_j}_{j=1}^{n}$ be a sequence of $c$-dimensional subspaces of $\bbF^d$.
For each $j=1,\dotsc,n$, let $\bfPhi_j$ be the $d\times c$ synthesis operator of an orthonormal basis for $\calU_j$, and let $\bfP_j=\bfPhi_j^{}\bfPhi_j^*$ be the corresponding orthogonal projection operator.
Then
\begin{equation}
\label{equation.equi-isoclinic Welch bound}
\max_{j\neq j'}\norm{\bfPhi_j^*\bfPhi_{j'}^{}}_{2}^2
\geq\tfrac{nc-d}{d(n-1)},
\end{equation}
where equality holds if and only if $\set{\calU_j}_{j=1}^{n}$ is an equi-isoclinic tight fusion frame (EITFF) for $\bbF^d$, namely if and only if both
\begin{enumerate}
\romani
\item
$\set{\calU_j}_{j=1}^{n}$ is a tight fusion frame,
namely there exists some $\alpha>0$ such that
$\displaystyle\sum_{j=1}^{n}\bfPhi_j^{}\bfPhi_j^{*}=\sum_{j=1}^{n}\bfP_j=\alpha\bfI$;
\item
$\set{\calU_j}_{j=1}^{n}$ is equi-isoclinic, namely there exists some $\sigma^2\geq0$ such that
$\bfPhi_{j'}^*\bfPhi_j^{}\bfPhi_j^*\bfPhi_{j'}^{}=\sigma^2\bfI$ for all $j\neq j'$, or equivalently, that $\bfP_{j'}\bfP_{j}\bfP_{j'}=\sigma^2\bfP_{j'}$ for all $j\neq j'$.
\end{enumerate}
In this case, $\set{\calU_j}_{j=1}^{n}$ is necessarily an ECTFF for $\bbF^d$, see Theorem~\ref{theorem.equi-chordal}, and $\alpha$ and $\sigma^2$ are necessarily
$\frac{nc}{d}$ and $\tfrac{nc-d}{d(n-1)}$, respectively.
\end{theorem}

In the special case where $c=1$, ECTFFs are equivalent to EITFFs since the induced $2$-norm of a $1\times 1$ matrix equals its Frobenius norm.
In fact, in this case, both ECTFFs and EITFFs correspond to ETFs:
letting $\set{\bfphi_j}_{j=1}^{n}$ be unit vectors chosen from the $1$-dimensional subspaces $\set{\calU_j}_{j=1}^{n}$ we have $\bfP_j=\bfphi_j^{}\bfphi_j^{*}$ and
$\Tr(\bfP_j\bfP_{j'})=\norm{\bfphi_j^*\bfphi_{j'}^{}}_{\Fro}^2=\abs{\ip{\bfphi_{j}}{\bfphi_{j'}}}^2$.

\section{Principal angles}

We have seen that ECTFFs give optimal packings of subspaces with respect to chordal distance.
Less obvious is what distance (if any) EITFFs are optimal packings with respect to.
To understand this better, we now discuss \textit{principal angles}.

Here, as before, let $\set{\bfPhi_j}_{j=1}^{n}$ be a sequence of $d\times c$ synthesis operators for orthonormal bases of a given sequence $\set{\calU_j}_{j=1}^{n}$ of $c$-dimensional subspaces of $\bbF^d$.
For each $j$ we have $\bfPhi_j^*\bfPhi_j^{}=\bfI$ and so $\norm{\bfPhi_j}_2\leq 1$.
As such, for any $j\neq j'$,
\begin{equation*}
\norm{\bfPhi_j^*\bfPhi_{j'}^{}}_2
\leq\norm{\bfPhi_j}_2\norm{\bfPhi_{j'}^{}}_2
\leq1.
\end{equation*}
This means that for any $j\neq j'$,
the singular values $\set{\sigma_{j,j',k}}_{k=1}^{c}$ of $\bfPhi_j^*\bfPhi_{j'}^{}$ are at most $1$.
As such, for any $j\neq j'$, there exists an increasing (non-decreasing) sequence of angles $\set{\theta_{j,j',k}}_{k=1}^{c}$ in $[0,\frac{\pi}2]$ such that $\sigma_{j,j',k}=\cos(\theta_{j,j',k})$ for all $k=1,\dotsc,c$.
For any $j\neq j'$, these angles are known as the \textit{principal angles} between $\calU_j$ and $\calU_{j'}$.
These angles are invariant with respect to the particular choice of orthonormal bases for $\set{\calU_j}_{j=1}^{n}$, since changing bases is equivalent to right-multiplying $\set{\bfPhi_j}_{j=1}^{n}$ by $c\times c$ unitary matrices,
which only affects the unitary terms in the singular value decomposition of $\bfPhi_j^*\bfPhi_{j'}^{}$, not its singular values.
One can express the chordal distance between any two subspaces in terms of their principal angles.
In particular, for any $j\neq j'$,
\begin{equation*}
\dist_{\rmc}^2(\calU_j,\calU_{j'})
=\tfrac12\norm{\bfP_j-\bfP_{j'}}_{\Fro}^2
=c-\norm{\bfPhi_j^*\bfPhi_{j'}^{}}_{\Fro}^2
=c-\sum_{k=1}^{c}\cos^2(\theta_{j,j',k})
=\sum_{k=1}^{c}\sin^2(\theta_{j,j',k}).
\end{equation*}
This ``chordal" notion of the distance between two subspaces has an advantage over other, more classical notions of distance, such as the \textit{geodesic distance} $(\sum_{k=1}^{c}\theta_{j,j',k}^2)^{\frac12}$:
ECTFFs give optimal packings with respect to the chordal distance,
whereas for any $c>1$, we are not aware of any practically-verifiable conditions that suffice to guarantee a given arrangement of subspaces is optimal with respect to geodesic distance.

Whereas ECTFFs are optimal packings with respect to the chordal distance,
EITFFs instead achieve equality in~\eqref{equation.equi-isoclinic Welch bound}:
\begin{equation*}
\tfrac{nc-d}{d(n-1)}
\leq\max_{j\neq j'}\norm{\bfPhi_j^*\bfPhi_{j'}^{}}^2
=\max_{j\neq j'}\max\set{\cos^2(\theta_{j,j',k})}_{k=1}^{c}
=\max_{j\neq j'}\cos^2(\theta_{j,j',1})
=1-\min_{j\neq j'}\sin^2(\theta_{j,j',1}).
\end{equation*}
That is, they ensure that $\min_{j\neq j'}\sin^2(\theta_{j,j',1})$ is as large as possible,
where for any $j\neq j'$, $\theta_{j,j',1}$ is the smallest principal angle between $\calU_j$ and $\calU_{j'}$;
Dhillon, Heath, Strohmer and Tropp refer to $\sin(\theta_{j,j',1})$ as the \textit{spectral distance} between $\calU_j$ and $\calU_{j'}$~\cite{DhillonHST08}.
Put more simply, EITFFs maximize the minimum principal angle of $\set{\calU_j}_{j=1}^{n}$, where this minimum is taken over all pairs of subspaces as well as all principal angles between them.

When viewed from the perspective of principal angles, EITFFs seem extremely special:
we need a collection of subspaces with such a high degree of symmetry that their orthogonal projection operators sum to a scalar multiple of the identity and such that any principal angle between any pair of subspaces is equal to any (other) principal angle between any (other) pair of subspaces.
EITFFs are so special, in fact, that one may reasonably doubt that nontrivial examples of them exist.
Nevertheless, they do:
if $\set{\bfphi_j}_{j=1}^{n}$ is any ETF for $\bbF^e$ and $c$ is any positive integer,
then letting $\bfI$ be the $c\times c$ identity matrix,
and letting $\bfPhi_j=\bfphi_j\otimes\bfI$ for all $j=1,\dotsc,n$,
we have that $\set{\calU_j}_{j=1}^{n}$, $\calU_j:=\operatorname{range}(\bfPhi_j)$, is an EITFF for $\bbF^{d}$ where $d=ce$.
To see this, note that letting $\bfPhi$ denote the $e\times n$ synthesis operator of the ETF $\set{\bfphi_j}_{j=1}^{n}$, the fusion frame operator of $\set{\bfPhi_j}_{j=1}^{n}$ is
\begin{equation*}
(\bfPhi\otimes\bfI)(\bfPhi\otimes\bfI)^*
=(\bfPhi\otimes\bfI)(\bfPhi^*\otimes\bfI)
=\bfPhi\bfPhi^*\otimes\bfI
=\tfrac ne(\bfI\otimes\bfI)
=\tfrac{nc}{d}\bfI,
\end{equation*}
while its fusion Gram matrix is
$(\bfPhi\otimes\bfI)^*(\bfPhi\otimes\bfI)=\bfPhi^*\bfPhi\otimes\bfI$,
meaning that for any $j\neq j'$,
\begin{equation*}
\norm{\bfPhi_j^*\bfPhi_{j'}^{}}_2^2
=\norm{(\bfphi_j\otimes\bfI)^*(\bfphi_{j'}\otimes\bfI)}_2^2
=\norm{\ip{\bfphi_j}{\bfphi_{j'}}\otimes\bfI}_2^2
=\tfrac{n-e}{e(n-1)}
=\tfrac{nc-ce}{ce(n-1)}
=\tfrac{nc-d}{d(n-1)}.
\end{equation*}
This trick allows one to use the growing list~\cite{FickusM15} of known ETF constructions to produce (infinite families of) nontrivial EITFFs,
all of which are also ECTFFs.
To our knowledge,
it is an open question whether every EITFF of $c$-dimensional subspaces of $\bbF^d$ has $n$ and $d$ parameters of the form $d=ce$ where there exists an $n$-vector ETF for $\bbF^e$.
More generally, it is an open question whether the dimension $c$ of the subspaces in an EITFF for $\bbF^d$ necessarily divides the ambient dimension $d$.

\section*{Acknowledgments}

This work was partially supported by NSF DMS 1321779, AFOSR F4FGA05076J002 and an AFOSR Young Investigator Research Program award.
The views expressed in this article are those of the authors and do not reflect the official policy or position of the United States Air Force, Department of Defense, or the U.S.~Government.

\end{document}